# NON-REPRODUCTIVE AND REPRODUCTIVE SOLUTIONS OF SOME MATRIX EQUATIONS

B.J. Malešević*, B.M. Radičić

* Faculty of Electrical Engineering, University of Belgrade, Serbia

*In this paper we analyzed solutions of some complex matrix equations related to pseudoinverses using the concept of reproductivity . Especially for matrix equation AXB = C it is shown that Penrose's general solution is actually the case of the reproductive solution.*

## Reproductive equations

A concept of reproductive equations was introduced by S.B. Prešić [4.] in 1968. The implementation of the concept of reproducibility was also considered by J.D. Kečkić ([6.], [7.]) and S.B. Prešić and J.D. Kečkić [12.].

**Definition.** *The reproductive equations* are the equations of the following form:

$$x = f(x)$$

where $x$ is a unknown, $S$ is a given set and $f : S \to S$ is a given function which satisfies

$$f \circ f = f\,. \qquad (1)$$

The condition (1) is called the *condition of reproductivity* because the function $f$ which satisfies the condition (1), after iteration, reproduces itself.

In literature, functions which satisfy the condition (1) are also called projectors or idempotent maps.

The following two statements (S.B. Prešić [4.]) give the fundamental properties of the reproductive equations.

**Theorem 1.** For any of the consistent equations $J(x)$ there is an equation of the $x = f(x)$ form, which is equivalent to $J(x)$ being in the same time reproductive as well. ◊

**Theorem 2.** If a certain equation $J(x)$ is equivalent to the reproductive one $x = f(x)$, then the general solution is given by the formula $x = f(y)$, for any of the values $y \in S$. ◊

The first implementation of reproducibility on some matrix equations was considered by S.B. Prešić ([2.], [4.]).

## Solutions of the matrix equation *AXB=C*

Let $m, n \in N$ and $C$ is the field of complex number. The set of all matrices of order $m \times n$ over $C$ we denote by $C^{m\times n}$. By $C_a^{m\times n}$ we denote the set of all matrices of order $m \times n$ over $C$ which are of rank *a*. Instead of $C^{m\times 1}$ we use the denotement $C_{\downarrow}^{m}$ while instead of $C^{1\times n}$ we use the denotement $C_{\rightarrow}^{n}$.

Let $A \in C^{m\times n}$**.** The solution of matrix equation

$$AXA = A \tag{2}$$

is called *{1}-inverse* of $A$ and it is denoted by $A^{(1)}$. The set of all {1}-inverses of $A$ is denoted by $A\{1\}$. More information about {1}-inverses and other types of inverses can be found in A. Ben-Israel and T.N.E. Greville [13.] and S.L. Campbell and C.D. Meyer [16].

Let $A \in C^{m\times n}$**,** $B \in C^{p\times g}$ and $C \in C^{m\times q}$**.** The following statement (R. Penrose [1.]) gives the necessary and sufficient condition for a matrix equation

$$AXB = C \tag{3}$$

is consistent. The formula of general solution of the matrix equation (3) is also given in the following statement.

**Theorem 3.** The matrix equation (3) is consistent iff for some of {1}-inverses $A^{(1)}$ and $B^{(1)}$ of the matrices $A$ and $B$ the following condition is true

$$AA^{(1)}CB^{(1)}B = C \tag{4}$$

The general solution of the matrix equation (3) is given by the formula

$$X = f(Y) = A^{(1)}CB^{(1)} + Y - A^{(1)}AYBB^{(1)} \tag{5}$$

where $Y$ is an arbitrary matrix corresponding dimensions. ◊

If the condition (4) is true for at least one choice of {1}-inverses $A^{(1)}, B^{(1)}$ of the matrices $A, B$ respectively, then the condition (4) is true for any choice of {1}-inverses $A^{(1)}, B^{(1)}$. The following statement is an extension of Theorem 3.

**Theorem 3'.** If the matrix $X_0 \in C^{n\times p}$ is any of the particular solutions of the matrix equation (3), then the general solution of the matrix equation (3) is given by the formula

$$X = g(Y) = X_0 + Y - A^{(1)}AYBB^{(1)} \tag{6}$$

for any choice of {1}-inverses $A^{(1)}, B^{(1)}$ of the matrices $A, B$ respectively and where $Y$ is an arbitrary matrix corresponding dimensions. The general solution (6) of the matrix equation (3) is reproductive one iff $X_0 = A^{(1)}CB^{(1)}$.

**Proof.** It is easily to show that the solution of the matrix equation is given by (6). On the contrary, let $X$ is any of the solutions of the matrix equation $(3)$. Based on Theorem 3, there is $Y$ so that $X = f(Y) = A^{(1)}CB^{(1)} + Y - A^{(1)}AYBB^{(1)}$. Then, for $Z = Y - X_0 + A^{(1)}CB^{(1)}$ it is true $X = f(Y) = f(Z + X_0 - A^{(1)}CB^{(1)})$ i.e.

$$X = g(Z) = f(Z + X_0 - A^{(1)}CB^{(1)}) = X_0 + Z - A^{(1)}AZBB^{(1)}.$$

From that we conclude that $(6)$ is the general solution of the matrix equation $(3)$.

Considering that the equality $g^2(Y) = g(Y) + (X_0 - A^{(1)}CB^{(1)})$ is valid, the formula of general solution (6) of the matrix equation (3) is reproductive iff $X_0 = A^{(1)}CB^{(1)}$. ◊

From the previous theorem we can see that Penrose's general solution (5) is the reproductive solution. In this paper we examine whether there are particular solutions $X_0$ which can not be represented in the form $X_0 = A^{(1)}CB^{(1)}$ (see the example at the end of this paper). The first appearances of these solutions were in $[2.]$. In that paper S.B. Prešić analyzed the matrix equation (2) in the case when $A$ is a square matrix and proved the following statement.

**Theorem 4.** For the square matrix $A \in C^{n\times n}$ and matrix $B \in A\{1\}$ the following is true:
1) $AX = 0 \Leftrightarrow \left(\exists Y \in C^{n\times n}\right)\ X = Y - BAY$, 2) $AX = A \Leftrightarrow \left(\exists Y \in C^{n\times n}\right)\ X = I + Y - BAY$,
3) $XA = 0 \Leftrightarrow \left(\exists Y \in C^{n\times n}\right)\ X = Y - YAB,$ 4) $XA = A \Leftrightarrow \left(\exists Y \in C^{n\times n}\right)\ X = I + Y - YAB,$
5) $AXA = A \Leftrightarrow \left(\exists Y \in C^{n\times n}\right)\ X = B + Y - BAYAB$. ◊

The general solutions 2), 4) and 5) of the previous theorem are non-reproductive and as such they do not directly appear according to Penrose's theorem. Let us notice that the formula (6) gives both Prešić's solutions of the given equations and Penrose's solutions of the given equations. In the paper $[8.]$ M. Haverić showed that we can get Penrose's solution from Prešić's solution. Namely, she proved the following statement

**Theorem 5.** For the square matrix $A \in C^{n\times n}$ and matrix $B \in A\{1\}$ the following is true:
1) $AX = 0 \Leftrightarrow \left(\exists Y \in C^{n\times n}\right)\ X = Y - BAY$, 2) $AX = A \Leftrightarrow \left(\exists Y \in C^{n\times n}\right)\ X = BA + Y - BAY$,
3) $XA = 0 \Leftrightarrow \left(\exists Y \in C^{n\times n}\right)\ X = Y - YAB,$ 4) $XA = A \Leftrightarrow \left(\exists Y \in C^{n\times n}\right)\ X = AB + Y - YAB,$
5) $AXA = A \Leftrightarrow \left(\exists Y \in C^{n\times n}\right)\ X = BAB + Y - BAYAB$. ◊

The matrix equation (3) is the subject of the contemporary research (see [10.],[11.],[14.],[15.]).

**Solutions of non-homogeneous linear systems**

In this part of the text we are going to show that the general solution of the consistent non-homogeneous linear system can be obtain by general {1}-inverse. Namely, for each matrix $A \in C^{m\times n}$ there are regular matrices $P \in C^{n\times n}$ and $Q \in C^{m\times m}$ such that

$$QAP = E_a = \begin{bmatrix} I_a & 0 \\ 0 & 0 \end{bmatrix}, \tag{7}$$

where $a = rank(A)$. Then, each {1}-inverse $A^{(1)}$ of the matrix $A$ can be represented in the following form

$$A^{(1)} = P \cdot \begin{bmatrix} I_a & U \\ V & W \end{bmatrix} \cdot Q \tag{8}$$

where $U = [u_{i,j}]$, $V = [v_{i,j}]$ and $W = [w_{i,j}]$ are matrices with total $k = m \cdot n - a^2$ mutually independent variables. In the paper $[5.]$ V. Perić mentioned that the representation (8) of the general {1}-inverse was given by C. Rohde in his doctoral thesis $[3.]$.

**Theorem 6.** Let $A \in C_a^{m\times n}$, $c \in C_{\downarrow}^{m}$ and

$$Ax = c \tag{9}$$

is a consistent non-homogeneous linear system. Then, there is {1}-inverse $A^{(1)}$ of the matrix $A$ in the form (8) with $n-a$ arbitrary variables such that a vector

$$x = A^{(1)}c \tag{10}$$

determines the general solution of linear system (9).

**Proof.** Let for the matrix $A$ the regular matrices $Q$ and $P$ be determined so that (7) is true. Let us also notice that the vector (10) always determines one solution of the linear system (9) for any choice of the variables $u_{i,j}, v_{i,j}$ and $w_{i,j}$ in the general {1}-inverse $A^{(1)}$. We conclude that

$$a = rank(A) = rank(Q\cdot A) \ \text{ and } \ a = rank([A|c]) = rank([Q\cdot A|Q\cdot c]) \tag{11}$$

using the consistency of the linear system (9). From (11) we can conclude that the vector $c' = Q\cdot c$ has zeroes on the last $m-a^{th}$ coordinates and there is $j^{th}$ coordinate $c_j'$ ($1\le j\le a$) such that $c_j' \ne 0$ because the linear system (9) is non-homogeneous. Therefore, in the case of the general inverse (8) the following is true

$$x = P\cdot\begin{bmatrix} I_a & U \\ V & W\end{bmatrix}\cdot Q\cdot c = P\cdot\begin{bmatrix} 1 & 0 & \cdots & 0 & u_{1,1} & \cdots & u_{1,m-a} \\ 0 & 1 & \cdots & 0 & u_{2,1} & \cdots & u_{2,m-a} \\ \vdots & \vdots & \ddots & \vdots & \vdots & \ddots & \vdots \\ 0 & 0 & \cdots & 1 & u_{a,1} & \cdots & u_{a,m-a} \\ v_{1,1} & v_{1,2} & \cdots & v_{1,a} & w_{1,1} & \cdots & w_{1,m-a} \\ \vdots & \vdots & \ddots & \vdots & \vdots & \ddots & \vdots \\ v_{n-a,1} & v_{n-a,2} & \cdots & v_{n-a,a} & w_{n.-a,1} & \cdots & w_{n-a,m-a}\end{bmatrix}\cdot\begin{bmatrix} c_1' \\ c_2' \\ \vdots \\ c_a' \\ 0 \\ \vdots \\ 0\end{bmatrix}$$

$$= P\cdot\begin{bmatrix} c_1' \\ c_2' \\ \vdots \\ c_a' \\ \sum_{i=1}^{a} c_i' v_{1,i} \\ \vdots \\ \sum_{i=1}^{a} c_i' v_{n-a,i}\end{bmatrix} = P\cdot\begin{bmatrix} c_1' \\ c_2' \\ \vdots \\ c_a' \\ 0 \\ \vdots \\ 0\end{bmatrix} + P\cdot\begin{bmatrix} 0 \\ 0 \\ \vdots \\ 0 \\ 1 \\ \vdots \\ 0\end{bmatrix}\cdot\tau_1 + \cdots + P\cdot\begin{bmatrix} 0 \\ 0 \\ \vdots \\ 0 \\ 0 \\ \vdots \\ 1\end{bmatrix}\cdot\tau_{n-a},$$

for $\tau_1 = \sum_{i=1}^{a} c_i' v_{1,i}$, …, $\tau_{n-a} = \sum_{i=1}^{a} c_i' v_{n-a,i}$. By choosing the submatrix $V$ in the form

$$V = \begin{bmatrix} 0 & \cdots & v_{1,j}/c_j' & \cdots & 0 \\ \vdots & & \vdots & & \vdots \\ 0 & \cdots & v_{n-a,j}/c_j' & \cdots & 0\end{bmatrix} \tag{12}$$

we get $\tau_1 = v_{1,j}$, … , $\tau_{n-a} = v_{n-a,j}$. From this we conclude that the directrix of the previously determined affine space $x = A^{(1)}c$ is of $n-a$ dimension. The vector $x = A^{(1)}c$, with $n-a$ arbitrary variables, determines the general solution of linear system (9). ◊

**Theorem 7.** Let $B \in C_b^{m\times n}$, $c \in C_{\rightarrow}^{n}$ and

$$xB = c \tag{13}$$

is a consistent non-homogeneous linear system. Then, there is {1}-inverse $B^{(1)}$ of the matrix $B$ in the form (8) with $m-b$ arbitrary variables such that a vector

$$x = cB^{(1)} \tag{14}$$

determines the general solution of linear system (13).

**Proof.** The proof is analogous to the previous proof. ◊

Using *the Kronecker* (*tensor*) *product of matrices* we can rewrite the matrix equation $AXB = C$ **(**A. Ben-Israel and T.N.E. Greville [13.], A.K. Jain [9.]) as

$$\left(A \otimes B^T\right) \; vec(X) = vec(C), \tag{15}$$

where $vec(X)$ denotes the vector of the matrix $X$ which is formed by writing the rows of the matrix $X$ into a single column vector. Namely, $vec$ determines the operator $vec_{m,n} : C^{m\times n} \rightarrow C_{\downarrow}^{mn}$ that is defined as follows $vec_{m,n}(x_{i,j}) = x_{(i-1)\cdot n+j}$. The inverse operator $mat_{m,n} : C_{\downarrow}^{mn} \rightarrow C^{m\times n}$ is defined as follows $mat_{m,n}(x_i) = x_{p,q}$, for $p = [i/n]+1$ and $q = i(\bmod n)$.

**Example.** Let be given the matrix equation:

$$AXB = C \quad \Leftrightarrow \quad \begin{bmatrix} 1 & 2 & 1 \\ 0 & 1 & 0 \\ 1 & 1 & 1 \end{bmatrix} X \begin{bmatrix} 1 & 1 \\ 1 & 1 \\ 2 & 2 \end{bmatrix} = \begin{bmatrix} -3 & -3 \\ -1 & -1 \\ -2 & -2 \end{bmatrix}. \tag{16}$$

Using the Kronecker product the matrix equation (16) may be considered in the form of equivalent linear system:

$$\left(A \otimes B^T\right) \; vec(X) = vec(C) \quad \Leftrightarrow \quad \begin{bmatrix} 1 & 1 & 2 & 2 & 2 & 4 & 1 & 1 & 2 \\ 1 & 1 & 2 & 2 & 2 & 4 & 1 & 1 & 2 \\ 0 & 0 & 0 & 1 & 1 & 2 & 0 & 0 & 0 \\ 0 & 0 & 0 & 1 & 1 & 2 & 0 & 0 & 0 \\ 1 & 1 & 2 & 1 & 1 & 2 & 1 & 1 & 2 \\ 1 & 1 & 2 & 1 & 1 & 2 & 1 & 1 & 2 \end{bmatrix} \cdot \begin{bmatrix} x_{1,1} \\ x_{1,2} \\ x_{1,3} \\ x_{2,1} \\ x_{2,2} \\ x_{2,3} \\ x_{3,1} \\ x_{3,2} \\ x_{3,3} \end{bmatrix} = \begin{bmatrix} -3 \\ -3 \\ -1 \\ -1 \\ -2 \\ -2 \end{bmatrix}. \tag{17}$$

Based on (8) the general {1}-inverses of $A=\begin{bmatrix}1&2&1\\0&1&0\\1&1&1\end{bmatrix}$ and $B=\begin{bmatrix}1&1\\1&1\\2&2\end{bmatrix}$ are given by the following matrices:

$$A^{(1)}=\begin{bmatrix}1-a+2b-c+e & -2+a-2b-d-e & a-2b-e\\ -b & 1+b & b\\ c-e & d+e & e\end{bmatrix},\ \left(a,b,c,d,e\in C\right)$$

and

$$B^{(1)}=\begin{bmatrix}1-g-2h-p+q+2r & g-q & h-r\\ p-q-2r & q & r\end{bmatrix},\ \left(g,h,p,q,r\in C\right).$$

The matrix

$$X_0=A^{(1)}CB^{(1)}=\begin{bmatrix}-1+3c+d+g+2h-3cg-6ch-dg-2dh & -g+3cg+dg & -h+3ch+dh\\ -1+g+2h & -g & -h\\ -3c+3cg+6ch-d+dg+2dh & -3cg-dg & -3ch-dh\end{bmatrix}$$

is one solution of the matrix equation (16) but we are going to show that it is not the general solution, $\left(c,d,g,h\in C\right)$.

Now, we will solve the linear system (17) using Theorem 6. Applying the elementary operations we get the regular matrices $Q\in C^{6\times 6}$ and $P\in C^{9\times 9}$ such that

$$Q\cdot\left(A\otimes B^T\right)\cdot P=E_r=\begin{bmatrix}I_r & 0\\ 0 & 0\end{bmatrix},$$

where $r=rank\left(A\otimes B^T\right)=rank\left(A\right)\cdot rank\left(B\right)$. Namely, we get that $r=2$ and

$$Q=\begin{bmatrix}1&0&0&0&0&0\\0&0&1&0&0&0\\-1&1&0&0&0&0\\0&0&-1&1&0&0\\-1&0&1&0&1&0\\-1&0&1&0&0&1\end{bmatrix},\ P=\begin{bmatrix}1&-2&-1&-2&0&0&-1&-1&-2\\0&0&1&0&0&0&0&0&0\\0&0&0&1&0&0&0&0&0\\0&1&0&0&-1&-2&0&0&0\\0&0&0&0&1&0&0&0&0\\0&0&0&0&0&1&0&0&0\\0&0&0&0&0&0&-1&0&0\\0&0&0&0&0&0&0&1&0\\0&0&0&0&0&0&0&0&1\end{bmatrix}.$$

Then,

$$C' = Q \cdot vec(C) = \begin{bmatrix} -3 & -1 & 0 & 0 & 0 & 0 \end{bmatrix}^T$$

and

$$vec(X) = (A \otimes B^T)^{(1)} vec(C) = P \begin{bmatrix} 1 & 0 & u_{1,1} & u_{1,2} & u_{1,3} & u_{1,4} \\ 0 & 1 & u_{2,1} & u_{2,2} & u_{2,3} & u_{2,4} \\ \frac{v_{1,1}}{-3} & 0 & w_{1,1} & w_{1,2} & w_{1,3} & w_{1,4} \\ \frac{v_{2,1}}{-3} & 0 & w_{2,1} & w_{2,2} & w_{2,3} & w_{2,4} \\ \vdots & \vdots & \vdots & \vdots & \vdots & \vdots \\ \frac{v_{7,1}}{-3} & 0 & w_{7,1} & w_{7,2} & w_{7,3} & w_{7,4} \end{bmatrix} C' = \begin{bmatrix} -1 - v_{1,1} - 2v_{2,1} - v_{5,1} - v_{6,1} - 2v_{7,1} \\ v_{1,1} \\ v_{2,1} \\ -1 - v_{3,1} - 2v_{4,1} \\ v_{3,1} \\ v_{4,1} \\ v_{5,1} \\ v_{6,1} \\ v_{7,1} \end{bmatrix}. \quad (18)$$

From (18) we get the matrix of the general solution

$$X = mat(vec(X)) = \begin{bmatrix} -1 - v_{1,1} - 2v_{2,1} - v_{5,1} - v_{6,1} - 2v_{7,1} & v_{1,1} & v_{2,1} \\ -1 - v_{3,1} - 2v_{4,1} & v_{3,1} & v_{4,1} \\ v_{5,1} & v_{6,1} & v_{7,1} \end{bmatrix}.$$

For $v_{1,1} = v_{2,1} = v_{6,1} = v_{7,1} = 1$ and $v_{3,1} = v_{4,1} = v_{5,1} = 0$, the matrix

$$X_1 = \begin{bmatrix} -7 & 1 & 1 \\ -1 & 0 & 0 \\ 0 & 1 & 1 \end{bmatrix}$$

is the solution of (16), but $X_1 \neq A^{(1)}CB^{(1)}$ for any choice of {1}-inverses $A^{(1)}, B^{(1)}$ because from $X_0 = X_1$ we obtain the contradiction $(g = 0$ and $(0 =) - 3cg - dg = 1)$. Finally, we proved that $X_1 \neq A^{(1)}CB^{(1)}$ for any choice of {1}-inverses $A^{(1)}$, $B^{(1)}$. ◊

**Acknowledgements**

The research is partially supported by the Ministry of Education and Science, Serbia, Grant No. 174032.